\documentclass[11pt]{amsart}
\usepackage{amssymb}

\newtheorem {theorem}{Theorem}
\newtheorem{lemma}{Lemma}

\newtheorem{definition}{Definition}

\newcommand{\inv}{^{-1}}
\newcommand{\abs}[1]{|#1|}
\newcommand{\eps}{{\varepsilon}}

\begin{document}

\begin{abstract}
  Equations in free groups have become prominent recently in
  connection with the solution to the well known Tarski Conjecture.
  Results of Makanin and Rasborov show that solvability of systems of
  equations is decidable and there is a method for writing down in
  principle all solutions. However, no practical method is known; the
  best estimate for the complexity of the decision procedure is
  P-space.

  The special case of one variable equations in free groups has been
  open for a number of years, although it is known that the solution
  sets admit simple descriptions. We use cancellation arguments to
  give a short and direct proof of this result and also to give a
  practical polynomial time algorithm for finding solution sets. One
  variable equations are the only general subclass of equations in
  free groups for which such results are known.

  We improve on previous attempts to use cancellation
  arguments by employing a new method of reduction motivated by
  techniques from formal language theory. Our paper is self contained;
  we assume only knowedge of basic facts about free groups.
\end{abstract}

\title{Solving One-Variable Equations in Free Groups}\author{Dimitri Bormotov
\ Robert Gilman \ Alexei Myasnikov}

\maketitle

\section{Introduction}

A one variable equation $E(x) = 1$ of degree $d$ in a finitely
generated free group $F$ is an expression of the form
\begin{equation}\label{equation}
  u_0 x^{\eps_0} u_1 x^{\eps_1} \ldots u_{d - 1}x^{\eps_{d-1}} = 1
\end{equation}
composed of elements $u_i \in F$, integers $\eps_i = \pm 1$ and a
symbol $x$ not in $F$. A solution to (\ref{equation}) is an element $g
\in F$ such that substitution of $g$ for $x$ yields $1$ in $F$.

Lyndon~{\cite{RCL1}} was the first to study equations of this sort.
He characterized solution sets in terms of parametric words. The
parametric words involved were simplified by
Lorents~{\cite{Lo1,Lo2}} and Appel~{\cite{Ap}}. However, Lorents
announced his results without proof, and Appel's published proof has
a gap (see \cite{CR}). A complete proof has been provided
recently by Chiswell and Remeslennikov~{\cite{CR}}.

Chiswell and Remeslennikov's novel analysis involves
algebraic geometry (\cite{BMR}, \cite{MR1}.) First they describe
the isomorphism types of the coordinate groups of irreducible
one-variable equations over $F$, and then they deduce the structure of
the solution sets. The latter part is easy, but the former requires
sophisticated techniques involving ultrapowers and Lyndon length
functions. The key point is that coordinate groups of irreducible
equations over $F$ are subgroups of the ultrapower $F^I/D$ of $F$ over
a countable set $I$ with a non-principal ultrafilter $D$.

One can view the group $F$ as a subgroup of $F^I/D$ under
the canonical diagonal embedding. From this point of view the
coordinate groups are precisely the finitely generated
subgroups of $F^I/D$ containing $F$ i.e., the so-called $F$-subgroups. In
particular up to isomorphism the coordinate groups of irreducible
one-variable equations over $F$ are the subgroups of $F^I/D$
of the form $\langle F, g\rangle$, $g \in F^I/D$.

Investigation of such $F$-subgroups of $F^I/D$ is not easy and
involves a careful analysis of Lyndon functions. (It might be
interesting to see whether it is easier to use free actions on
$\Lambda$-trees.) The computations can be simplified by employing a
result from \cite{KM} which states that the coordinate groups of
irreducible varieties are precisely the finitely generated
$F$-subgroups of the free exponential Lyndon group
$F^{\mathbb{Z}[t]}$. As this group is the union of an infinite
ascending chain of extensions of centralizers of $F$ \cite{MR2}, one
can use Bass-Serre theory to study $F$-subgroups of
$F^{\mathbb{Z}[t]}$.

Chisewell and Remeslennikov's method is very powerful and
potentially useful for more than just free groups.  However, it does
have the disadvantage of not giving an algorithm for explicitly describing
the set of solutions.

This paper is a refinement and extension of~\cite{GiMy} where
results from formal language theory are used to describe solution
sets of one-variable equations in free groups. As it turns out,
formal language techniques are not required; straightforward
cancellation arguments suffice. It seems likely that these arguments
can be extended to other groups admitting suitable (not necessarily
Lyndon) length functions. The main advantage of this method is that
it is short and yields a polynomial time algorithm for producing a
description of all solutions. This algorithm has been implemented by
the first author \cite{Bo}.

\begin{theorem}\label{th:main}
  \label{theorem} The solution set for a one variable equation of positive
  degree in a free group $F$ is a finite union of sets $uv^i w$
  where $u, v, w \in F$ and $i$ ranges over all integers. There is a
  polynomial time algorithm for finding these sets.
\end{theorem}

Let $\Sigma$ be a set of free generators for $F$ together with their
inverses, and let $\Sigma^{\ast}$ be the free monoid over $\Sigma$.
We consider Equation (\ref{equation}) in terms of words in
$\Sigma^{\ast}$.  Each coefficient $u_i$ is represented by
a freely reduced word (also denoted $u_i$) in $\Sigma^{\ast}$. From
this point of view $E (x) = u_0 x^{\eps_1} u_2 x^{\eps_2} \ldots
u_{d - 1} x^{\eps_d} u_d$ is a word in the free monoid over $\Sigma
\cup \{x, x \inv \}$, and a solution to $E (x) = 1$ is a word $s \in
\Sigma^{\ast}$ such that $E (s)$ is freely equal to the empty word.
The first assertion of Theorem~\ref{th:main} is equivalent to saying
that for some finite union of sets of words $uv^i w$ the solutions set
consists of all words freely equal to elements of the finite union.
A set $uv^i w$ is called a parametric word.

We assume without loss of generality that $E (x)$ is freely reduced,
and call $d$ the degree of $E (x)$. If $d = 0$, then $E (x) = u_0$.
In this case the solution set is empty if $E (x) \ne 1$ and all of
$\Sigma^{\ast}$ if $E (x) = 1$.  If the equation has degree one, it is
easy to find its unique solution. From now on we consider only
equations of degree at least two.


We begin with some lemmas on cancellation, after which we find a finite
number of parametric words $uv^i w$ and $uv^i wr^j s$ which contain
all solutions to $E(x) = 1$ up to free equivalence. Next we show that
two parameters are
not required and that $uv^i w$ is either a solution for all integers
$i$ or for an effectively determined finite subset. At the end we
present the algorithm and estimate its time complexity.

To explain our argument in more detail we require a few definitions.
For any (word) $g \in F$ we say that the $i$th occurrence of $g$
 {\em cancels out} in $E(g)$
 if there exists a way to freely reduce $E(g)$ such
 that all letters from $g^{\epsilon_i}$ cancel out during this
 reduction process.

We say that $g$ is a {\em pseudo-solution} of $E(x) = 1$ if some
occurrence of $g$ cancels out in $E(g)$.  Obviously every solution
of $E(x) = 1$ is also a pseudo-solution of $E(x) = 1$. However, unlike
solutions, pseudo-solutions admit a nice reduction
theory.

Our key idea is to study pseudo-solutions of equations instead of
solutions. The first result in this direction (stated in \cite{GiMy}
in a slightly different form) reduces the situation to cubic
equations. Namely, Lemma \ref{reduction} shows that if $g$ is a
pseudo-solution of $E(x) = 1$ in $F$ then $g$ is a pseudo-solution of
a cubic equation of the type
 $$x^{\eps_{j - 1}} u_{j} x^{\eps_j} u_{j+1} x^{\eps_{j + 1}},$$ where
    $0< j < d$ and indices are read modulo $d$ (so $u_d=u_0$.)
Next in Lemma \ref{parametric words} we show that pseudo-solutions
of cubic equations are in fact pseudo-solutions of some particular
quadratic equations which one can find effectively.
Finally, Lemmas \ref{sus} and \ref{susinv}  give a precise
description of pseudo-solutions of quadratic one-variable equations
over $F$ in terms of parametric words. Combining all these results
we obtain description of all pseudo-solutions of $E(x) =
1$ in terms of parametric words in two parameters.

The rest of our proof explains precisely how to use only one parameter
to describe solutions of $E(x) = 1$. The method of big powers (see
\cite{BMR2}) is the key tool in the second part. This means that the
argument is rather general - it works in many other groups that
satisfy the big powers condition (see \cite{KvM}), for example
torsion-free hyperbolic groups.

  One-variable equations are the only general
  class of equations in free groups for which a good description of
  solution sets as well as a practical (polynomial time) algorithm are
  known. In his seminal paper \cite{Mak82} Makanin proved decidability
  of the Diophantine problem in free groups $F$ (whether or not a
  given equation has a solution in $F$); however, his original
  algorithm is very inefficient - not even primitive recursive (see
  \cite{KoKoPa}). In the fundamental paper \cite{Razborov} Razborov
  gave a  description of solution sets of arbitrary equations
  in $F$.  Though this description is extremely complicated, it was useful
  in the solution of several deep problems in group theory
  \cite{KM2,KM3,BKM} including the Tarski's problems \cite{KM4}. In
  another paper \cite{Razborov2} Razborov showed that, in general,
  there is no easy description of solutions sets of equations in
  $F$. Later, Plandowski gave a much improved $P$-space version of the
  decision algorithm for equations in free monoids \cite{Plandowski},
  and Gutierrez devised a 
  $P$-space algorithm for the decision problem for equations in free
  groups \cite{Gu}.  Recent results \cite{DGH} due to Diekert,
  Gutierrez, and Hagenah, indicate that the decision problem for
  equations in free groups might be $P$-space-complete, though nothing
  definite has been proven so far. These results on the complexity of
  the decision problem for equations in free groups and for their
  solution sets make the existence of 
  subclasses of equations admitting polynomial decision algorithms and 
  descriptions of solutions sets in closed form, all
  the more remarkable.


\section{Cancellation Lemmas}

As above $\Sigma$ is a set of free generators and their inverses for a free
group $F$, and $\Sigma^{\ast}$ is the free monoid over $\Sigma$. Let $p, q, r,
s, t, u, v, w$ be words in $\Sigma^{\ast}$. We write $u \sim v$ if $u$ is
freely equal to $v$, and $u \to v$ if $u$ can be reduced to $v$ by
cancellation of subwords $aa \inv$, $a \in \Sigma$. In particular $u \to u$.
The empty word is denoted $1$, and the length of $u$ is $\abs{u}$. Recall that
for any word $u$ there is a unique irreducible word $v$ such that $u \to v$,
and further $u \sim w$ if and only if $w \to v$.

We introduce some additional notation.

\begin{definition}\label{notation} Let $w$ be any word.
\begin{enumerate}
\item $w'$ stands for an arbitrary prefix of $w$ and $w''$ for an
  arbitrary suffix.
\item $\abs w_c$ is the length of a cyclicly reduced word conjugate to $w$.
\end{enumerate}
\end{definition}

\begin{lemma}
  \label{factorization}If $v \to u$ and $u = u_1 u_2 \cdots u_m$, then $v =
  v_1 v_2 \cdots v_m$ with $v_i \to u_i$.
\end{lemma}

\begin{proof}
  Use induction on $n$, the number of cancellations necessary to reduce $v$ to
  $u$. If $n = 0$, then $u = v$ and there is nothing to prove. Otherwise let
  the first reduction be $v \to w$. By induction $w = w_1 w_2 \cdots w_m$ with
  $w_i \sim u_i$. As $v$ is obtained from $w$ by inserting a subword $aa \inv$
  into some $w_i$ or appending it to the beginning or end of some $w_i$, $v$
  has the desired factorization.
\end{proof}

\begin{lemma}
  \label{cancellation}Consider a fixed sequence of cancellations which reduces
  $u$ to $v$. If two particular letters of $u$ cancel at some point in the
  sequence, then either they are adjacent in $u$ or the subword between them
  has been reduced to $1$ by previous cancellations.
\end{lemma}

\begin{proof}
  Use induction on the length of the cancellation sequence.
\end{proof}

Now we slightly generalize the definition of a
pseudo-solution of equation to the following situation.
\begin{definition}
  A subword $s$ of $w$ is a pseudosolution if there is a sequence of
  cancellations in $w$ which consumes all letters in $s$.
\end{definition}

We are dealing with words over $\Sigma$, not group elements. For
example $s = ab \inv$ is a pseudosolution of $asba \inv a$ but not of $asb$.
The next two lemmas can be proved by straightforward induction on the length
of an appropriate cancellation sequence.

\begin{lemma}
  \label{split}Suppose $s$ is a pseudosolution of $w = usv$, then $s = s_1
  s_2$ with $s_1$ a pseudosolution of $us_1$ and $s_2$ a pseudosolution of
  $s_2 v$.
\end{lemma}

\begin{lemma}
  \label{subword}Let $s$ be a pseudosolution of $w$, and fix a cancellation
  sequence. The smallest subword of $w$ which contains $s$ and all letters in
  $w$ canceling with letters of $s$ is freely equal to $1$.
\end{lemma}

\begin{lemma}
  \label{subword2}A subword $s$ of $w$ is a pseudosolution if and only if
  there is a word $t$ such that $s$ is a subword of $t$, $t$ is a subword of
  $w$, and $t \sim 1$.
\end{lemma}

\begin{proof}
  If $t$ exists, then $t \sim 1$ implies $t \to 1$ whence $t$ and all its
  subwords are pseudosolutions of $w$. For the converse apply
  Lemma~\ref{subword}.
\end{proof}

\begin{lemma}
  \label{sus}If $us$ and $sw$ are irreducible and if either occurrence of $s$
  is a pseudosolution of $usvsw$, then $s \sim v_3 \inv v_1 \inv$ for some
  factorization $v = v_1 v_2 v_3$.
\end{lemma}

\begin{proof}
  We argue by induction on $n$, the length of a cancellation sequence. If $n =
  0$, then $s = 1$ in which case we take $u_1 = u_3 = 1$ and $u_2 = u$. If $v
  = 1$, then $usvsw = ussw$. As $us$ and $sw$ are irreducible, the only
  reduction possible involves cancellation at the boundary between $us$ and
  $sw$. It follows that $ss \sim 1$, whence $s \sim 1$.

  Assume $n > 0$ and $v \ne 1$. If the first reduction is within $v$, then $v
  \to p$ and by induction $s \sim p_3 \inv p_1 \inv$ for some factorization $p
  = p_1 p_2 p_3$. Lemma~\ref{factorization} implies $v = v_1 v_2 v_3$ with
  $v_i \sim p_i$ and $s \sim v_3 \inv v_1 \inv$.

  The remaining possibilities are cancellation at the boundary between $s$ and
  $v$ or the boundary between $v$ and $s$. Consider the first case; the second
  is similar. We have $s = ta \inv$ and $v = ap$ for some letter $a$ and words
  $t$ and $p$. The induction hypothesis applied to $utpt (a \inv w)$ yields $p
  = p_1 p_2 p_3$ and $t \sim p_3 \inv p_1 \inv$. But then $v = ap = (ap_1) p_2
  p_3$ and $s = ta \inv \sim p_3 \inv (ap_1) \inv$ as desired.
\end{proof}

\begin{lemma}
  \label{susinv}Suppose $v \not\sim 1$. If $us$ and $s \inv w$ are
  irreducible and $s$ or $s \inv$ is a pseudosolution of $usvs \inv w$, then
  $s \sim v'' v^k$ for some integer $k$. (See
  Definition~\ref{notation}.) Likewise if $us \inv$ and $sw$ are
  irreducible and $s$ or $s \inv$ is a pseudosolution of $us \inv vsw$, then
  $s \sim v^k v'$.
\end{lemma}

\begin{proof}
  Consider the first part; as before use induction on $n$, the number of
  cancellations. If $n = 0$, then $s = 1$. Take $v_1 = v, v_2 = 1$ and $k =
  0$. Otherwise the first reduction is either within $v$ or at one end or the
  other of $v$. In the first case $v \to v'$, and the induction hypothesis
  applied to $usv' s \inv w$ yields the desired result.

  Suppose then that there is a reduction at the left end of $v$; the other
  case is similar. We have $s = ta \inv$, $v = ap$, and application of the
  induction hypothesis to $ut (pa) t \inv w$ yields $pa = p_1 p_2$ and $t \sim
  p_2 (pa)^k$. It follows that $s \sim p_2 (pa)^k a \inv \sim p_2 a \inv a
  (pa)^k a \inv \sim p_2 a \inv v^k$. If $p_2 \ne 1$, then $p_2 = v_2 a$ for
  some suffix $v_2$ of $v$ whence $s \sim v_2 v^k$. If $p_2 = 1$, then $s \sim
  a \inv v^k \sim a \inv vv^{k - 1} \sim pv^{k - 1}$. As $p$ is a suffix of
  $v$, the first assertion holds. The second assertion follows from the first
  upon replacement of $s$ by $s \inv$.
\end{proof}

\begin{lemma}
  \label{ts}If $s$ is a pseudosolution of $tus$, $t$ is a pseudosolution of
  $tvs$, and $st$ is irreducible, then $s \sim (v \inv u)^i (v \inv u)'$ and
  $t \sim (vu \inv)'' (vu \inv)^j$ for some integers $i, j$.
\end{lemma}

\begin{proof}
  Application of Lemma~\ref{subword} to $tus$ implies either $u = u_1 u_2$
  with $u_2 s \sim 1$ or $t = t_1 t_2$ with $t_2 us \sim 1$. Consider the
  first case. We have $s \sim u_2 \inv \sim (v \inv u) \inv (v \inv u_1)$ as
  required. Further $u_2 \inv \to s$ implies that $t$ is a pseudosolution of
  $tvu_2 \inv$ and hence of $tvu \inv$. Thus either $v = v_1 v_2$ with $t \sim
  v_1 \inv$ or $u = u_3 u_4$ with $t \sim (vu_4 \inv) \inv$. But then $t \sim
  (v_2 u \inv) (vu \inv) \inv$ or $t \sim (u_3 \inv) (vu \inv) \inv$, and we
  see that $t$ has the right form. A similar analysis starting with starting
  with $tvs$ also works.

  It remains to consider the case $t = t_1 t_2$ with $t_2 us \sim 1$ and $s =
  s_1 s_2$ with $tvs_1 \sim 1$. Suppose $u \sim v$. We have $t = t_1 t_2$ with
  $t_2 us \sim 1$ and $s = s_1 s_2$ with $tus_1 \sim 1$. If $t_1 = 1$, then
  $tus \sim 1$ implies $st \sim u \inv$. As $st$ is irreducible,
  Lemma~\ref{factorization} yields $s \sim (u \inv)' = (v \inv)'$ and $t \sim
  (u \inv)''$ which is included in $i = j = 0$. If $t_1 \ne 1$, it follows
  from $t_2 us_1 s_2 \sim 1 \sim t_1 t_2 us_1$ that $t_1 \sim s_2$. As $st$ is
  irreducible, $t_1$ and $s_2$ are too. Thus $t_1 = s_2 \ne 1$. Hence $s_1 s_2
  t_2 = s_1 t_1 t_2$ is irreducible. But then $s_1 s_2 t_2 \sim u \inv \sim
  s_1 t_1 t_2$ implies $s \sim (v \inv)', t \sim (u \inv)''$ as before.

  Finally suppose $t = t_1 t_2$ with $t_2 us \sim 1$, $s = s_1 s_2$ with
  $tvs_1 \sim 1$, and $u \not{\sim} v$. From $t_2 us \sim 1$ we deduce $u \inv
  t_2 \inv \to s$. Hence $t$ is a pseudosolution of $tvu \inv t_2 \inv$ and
  all the more of $tvu \inv t_2 \inv t_1 \inv = tvu \inv t \inv$. Likewise $s$
  is a pseudosolution of $s \inv v \inv us$. We are done by
  Lemma~\ref{susinv}.
\end{proof}

\begin{lemma}
  \label{stus} Let $st$ be irreducible. If the right-hand occurrence
  of $s$ is a
  pseudosolution in $stus$ but not in $tus$, then $st \sim u_3 \inv u_1 \inv$
  for some factorization $u = u_1 u_2 u_3$. Likewise if the left-hand
  occurrence of $t$ is a pseudosolution in $tvst$ but not in $tvs$, then $st
  \sim v_3 \inv v_1 \inv$ for some factorization $v = v_1 v_2 v_3$.
\end{lemma}

\begin{proof}
  Consider the first part; the second is treated similarly. We have $s = pq$
  with $q \ne 1$ and $qtus \sim 1$. Since $st$ is irreducible, so is $qt$. It
  follows that $t$ is a pseudosolution of $tus$. If $t$ is not a
  pseudosolution of $tu$, then $s = ef$ with $tue \sim 1$. But then $qf \sim
  1$ forces $f$ to be a pseudosolution of $pqf = sf = eff$, and
  Lemma~\ref{sus} implies $f \sim 1$. Consequently $tuef = tus \sim 1$
  contrary to our hypothesis that $s$ is not a pseudosolution of $tus$.

  It remains to deal with the possibility that $t$ is a pseudosolution of
  $tu$. In this case $u = u_1 u_2$ with $t \sim u_1 \inv$. It follows that
  $qu_2 s \sim 1$ whence the right-hand occurrence of $s$ is a pseudosolution
  in $su_2 s$. An application of Lemma~\ref{sus} completes the proof.
\end{proof}

\begin{lemma}\label{powers1} Suppose $p^iuq^j\sim v$, $i,j\ge 0$, and
  $i\abs p_c + j \abs q_c \ge 2\abs p + 2\abs q + \abs u +
  \abs v$. (Recall Definition~\ref{notation}.) Further assume that $p$
  and $q$ are not freely equal to proper powers. Under these
  conditions $uqu\inv\sim p\inv$.
\end{lemma}

\begin{proof} Assume the Lemma holds when both $p$ and $q$ are cyclicly
  reduced, and consider the case that they are not. Free reduction of
  $p$ and $q$ yields reduced words $p_1p_2p_1\inv\sim p$,
  $q_1q_2q_1\inv\sim q$ with $p_2, q_2$ cyclicly reduced. Hence
  $p_2^i(p_1\inv u q_1) q_2^j\sim p_1\inv vr_1$. Rewriting $p_1\inv u
  q_1$ as $u_2$ and $p_1\inv vq_1$ as $v_2$ we obtain $p_2^iu_2
  q_2^j\sim v_2$. As $i\abs{p_2}_c + j\abs{q_2}_c =i\abs
  p_c + j\abs q_c \ge  2\abs p + 2\abs q + \abs u +
  \abs v \ge 2\abs{p_2} + 2\abs{u_2} + \abs{q_2} + \abs{v_2}$, we have
  $u_2q_2u_2\inv \sim p_2^{\pm 1}$. Hence $uqu\inv\sim
  (p_1u_2q_1\inv)(q_1q_2q_1\inv)(q_1u_2\inv p_1\inv)\sim p_1u_2q_2u_2\inv
  p_1\inv \sim p_1p_2^{\pm 1}p_1\inv\sim p\inv$.

  It remains to deal with the case that $p$ and $q$ are cyclicly
  reduced. Without loss of generality assume that $u$ and $v$ are
  freely reduced and $i,j\ge 0$. Thus there is a sequence of
  $(1/2)(\abs{p^iuq^j}-\abs v) = (1/2)(i\abs p + j\abs q + \abs u -
  \abs v) \ge \abs p + \abs q + \abs u$ cancellations which reduces
  $p^iuq^j$ to $v$.

  Since cancellation can occur only at either end of $u$, the first
  $\abs u$ cancellations must consume $u$. In other words $u$ cancels
  with a suffix of $p^i$ and a prefix of $q^j$. For some
  factorizations $p=p_1p_2$ and $q=q_1q_2$ we have $u=(p_2p^{i_2})\inv
  (q^{j_1}q_1)\inv$ with $i=i_1+1 + i_2$ and
  $j=j_1+1+j_2$. Consequently $p^{i_1}p_1q_2q^{j_2}$ admits at least
  $\abs p + \abs q$ cancellations. Thus the infinite sequences $q_2
  q_1 q_2 q_1 \cdots$ and $ p_1\inv p_2\inv p_1 \inv p_2 \inv \cdots$
  have the same prefix of length $\abs{p} + \abs{q}$. As these
  sequences have periods $\abs{p}$ and $\abs{q}$ respectively, they
  are identical by~\cite[Theorem 1]{FW}. But then the fact that
  $(p_1\inv p_2\inv)^{\abs q}$ and $(q_2 q_1)^{\abs p}$ have the same
  length implies that they are equal. Since $p$ and $q$ are not
  proper powers, neither are  $(p_1\inv p_2\inv)$ and $q_2 q_1$. It
  follows that  $p_1\inv p_2\inv = q_2 q_1$, and this equation implies
  in a straightforward way that $uqu\inv\sim p\inv$.
\end{proof}

\begin{lemma}\label{powers2} Suppose that $q^j$ is a
  pseudosolution of $p^iuq^jvr^k$, $\abs j \abs q_c \ge 7(\abs p + \abs
  u + \abs q + \abs v + \abs r)$; and $p,q,r$ are not proper
  powers. Then either $q\sim 1$ or $\abs i\ge 1$ and $u\inv p u \sim
  q^{\pm 1}$ or $\abs k\ge 1$ and $vrv\inv\sim q^{\pm 1}$.
\end{lemma}

\begin{proof} Without loss of generality assume $i,j,k\ge 0$. By
  Lemma~\ref{split} $q$ factors as $q_1q_2$ and $q^j$
  factors as $(q^{j_1}q_1)(q_2q^{j_2})$ in such a way that
  $q^{j_1}q_1$ is a pseudosolution of $p^iuq^{j_1}q_1$, and
  $q_2q^{j_2}$ is a pseudosolution of $q_2q^{j_2}vr^k$. Clearly one of
  $j_1, j_2$ is no smaller than $(j-1)/2$. Assume it is $j_1$; the
  argument is similar in the other case.

  By Lemma~\ref{subword} $q^{j_1}q_1$ extends to a suffix of
  $p^iuq^{j_1}q_1$ which is freely equal to $1$. If that suffix is
  contained in $uq^{j_1}q_1$, then $u=u_1u_2$ with $q^{j_1}\sim
  u_2\inv q_1\inv$. Hence  $q^{j_1}$ freely reduces to a word $w$ with
  $\abs w \le \abs u + \abs q$. On the other hand $\abs w \ge
  j_1\abs q_c \ge .5(j-1)\abs q_c\ge 3.5(\abs p + \abs u + \abs
  q)- .5 \abs q \ge 3(\abs p + \abs q + \abs u)$. But then $\abs p =
  \abs q =\abs u =0$, which implies $q\sim 1$.

  It remains to consider the case that the suffix is not
  contained in $us^{j_1}s_1$. In particular $i\ge 1$.
  For some factorization $p=p_1p_2$ and $m\le i$ we have
  $p_2p^muq^{j_1}q_1\sim 1$. Thus $p^muq^{j_1}\sim p_2\inv
  q_1\inv$. As above $j_1\abs q_c\ge 3(\abs p + \abs q + \abs u) \ge
  2\abs p + 2\abs q + \abs u + \abs{p_2\inv q_1}$. Lemma~\ref{powers1}
  applies and yields $uqu\inv \sim p^{\pm 1}$.
\end{proof}


\section{Parametric Words}

In this section we show how to find a finite set of words and parametric words
$uv^i wr^j s$ which together contain all solutions to
Equation~(\ref{equation}).

Let $s$ be any freely reduced word which is a solution to
Equation~(\ref{equation}). Substitution of $s$ for $x$ yields a word
\begin{equation}
  \label{solution} E(s) = u_0 s^{\eps_0} \ldots u_{d - 1} s^{\eps_{d-1}}
\end{equation}
such that $E (s) \to 1$.

Fix a sequence of cancellations which reduces $E (s)$ to $1$, and let
$s^{\eps_j}$ be the first of the subwords $s^{\pm 1}$ to be consumed. If there
is a tie, pick either subword. Observe that the letters in $s^{\eps_j}$ must
cancel with nearby letters in $E (s)$. If a letter in $s^{\eps_j}$ canceled
to the right of $s^{\eps_{j + 1}}$, then by Lemma~\ref{cancellation}
$s^{\eps_{j + 1}}$ would disappear before $s^{\eps_j}$. Likewise no letter of
$s^{\eps_{j + 1}}$ cancels to the left of $s^{\eps_{j - 1}}$. We have the
following result.

\begin{lemma}\label{reduction} One of the following holds.
  \begin{enumerate}
    \item $s^{\eps_0}$ is a pseudosolution of $u_0 s^{\eps_0} u_1 s^{\eps_1}$;

    \item For some $j$ strictly between $0$ and $d-1$, $s^{\eps_j}$ is a
    pseudosolution of $s^{\eps_{j - 1}} u_{j} s^{\eps_j} u_{j+1} s^{\eps_{j +
    1}}$;

    \item $s^{\eps_{d-1}}$ is a pseudosolution of $s^{\eps_{d - 2}} u_{d-1}
    s^{\eps_{d-1}}$.
  \end{enumerate}
\end{lemma}

It is convenient to use the following immediate consequence of
Lemma~\ref{reduction}.

\begin{lemma}\label{cyclic-reduction} For some $j$ between $0$ and
  $d-1$, $s^{\eps_j}$ is a
    pseudosolution of $s^{\eps_{j - 1}} u_{j} s^{\eps_j} u_{j+1} s^{\eps_{j +
    1}}$. Here indices are read modulo $d$; e.g., $u_d=u_0$.
\end{lemma}

It follows from Lemma~\ref{cyclic-reduction} that application of the
following lemma to all successive pairs of coefficients $u = u_i$, $v
= u_{i + 1}$ (with indices read modulo $d$) yields a set of words and
parametric words containing $s$
or $s \inv$ for every solution $s$ to Equation~\ref{equation}.

\begin{lemma}\label{parametric words}
  If $\alpha, \beta = \pm 1$ and $s$ is an irreducible pseudosolution to
  $s^{\alpha} usvs^{\beta}$, then one of the following holds. (Recall
  Definition~\ref{notation}.)
  \begin{enumerate}
    \item \label{i}$s \sim (v \inv u)^i (v \inv u)' (vu \inv)'' (vu \inv)^j$;

    \item \label{ii}$s \sim (u \inv)' (u \inv)''$ or $(v \inv)' (v \inv)''$;

    \item \label{iii}$s \sim (u \inv)' v^i v' v'' v^j$ or $u^i u' u''
    u^j (v \inv)''$;

    \item \label{iv}$s \sim u^i u' v'' v^j$.
  \end{enumerate}
\end{lemma}

\begin{proof}
  By Lemma~\ref{split} $s = s_1 s_2$ with $s_1$ a pseudosolution of
  $s^{\alpha} us_1$ and $s_2$ a pseudosolution of $s_2 vs^{\beta}$. There are
  four cases. First if $\alpha = - 1, \beta = - 1$, Lemma~\ref{susinv} applied
  to $s_2 \inv s_1 \inv u \underline{s_1}$ and $\underline{s_2} vs_2 \inv s_1
  \inv$ yields~(\ref{iv}).

  If $\alpha = \beta = 1$, we have $s_1 s_2 u \underline{s_1}$ and
  $\underline{s_1} vs_1 s_2$ where the pseudosolutions are underlined. It may
  happen that $s_1$ is pseudosolution of $s_2 us_1$ and $s_2$ is a
  pseudosolution of $s_2 vs_1$. In this case Lemma~\ref{ts} applies
  and~(\ref{i}) holds. Otherwise either $s_1$ is not a pseudosolution of $s_2
  us_1$ or $s_2$ is not a pseudosolution of $s_2 vs_1$. In both cases
  Lemma~\ref{stus} implies~(\ref{ii}).

  Suppose $\alpha = 1, \beta = - 1$. In this case $s_1 s_2 u \underline{s_1}$
  and $\underline{s_2} vs_2 \inv s_1 \inv$. By Lemma~\ref{stus} either $s$ in
  included in~(\ref{ii}) or $s_2 u \underline{s_1}$ whence $s_1$ is freely
  equal to the inverse of a suffix of $s_2 u$. Equivalently $s_1$ is freely
  equal to a prefix of $(s_2 u) \inv$. But Lemma~\ref{susinv} implies $s_2
  \sim v_2 v^j$ for some integer $j$ and factorization $v = v_1 v_2$. It
  follows from Lemma~\ref{factorization} that $s_1$ is freely equal to a
  prefix of $(v_2 v^j u) \inv$. Consideration of the possible cases
  yields~(\ref{iii}).

  A similar argument works when $\alpha = - 1, \beta = 1$ and shows
  that~(\ref{ii}) or~(\ref{iii}) holds.
\end{proof}


\section{Solutions}

In order to find all solutions to Equation~(\ref{equation}) we need to test
the possibilities given by Lemma~\ref{parametric words}. It is
straightforward to test the single words; the parametric words require
more work. They have the form $rp^i sq^j t$. Without
loss of generality we assume that $p$ and $q$ are not proper
powers. By introducing words of the form $rp^is$ we may assume
$p\not\sim 1 \not\sim q$.

Consider $rp^is$. Substitute
$rys$ for $x$ in Equation~\ref{equation} to obtain an equation
$E'(y)=v_0y^{\eps_0}\cdots v_{d-1}y^{\eps_{d-1}}$ in the
indeterminate $y$ with coefficients $v_j$ of the form $su_jr$,
$su_js\inv$ etc. Note that $rp^is$ is a solution of $E(x)$ if and only
if $p^i$ is a solution of $E'(y)$. Also the sum of the lengths of the
coefficients of $E'(y)$ is $\abs{v_0\cdots v_{d-1}}=\abs{u_0\cdots
  u_{d-1}}+ d\abs{rs}$. Denote this number by $K_1$.

If a coefficient $v_{j}$ commutes with $p$, i.e.\ $v_jp\sim pv_j$,
then the subword $y^{\eps_{j-1}}v_jy^{\eps_j}$ of $E'(y)$ may be
replaced by $v_jy^{\eps_{j-1}+\eps_j}$ without affecting the set of
$i$'s for which $p^i$ is a solution. This is true even if indices are
read modulo $d$. The coefficients in $E'(y)$ will change, but
$E'(1)=v_0\cdots v_{d-1}$ remains constant. In particular the sum of
the length of the coefficients is still $K_1$.

Continue replacements of this sort until reaching an equation of the
form $E''(y) = w_0y^{k_0}\cdots w_my^{k_m}$ with $m$ minimal. It may
be that $m=0$ and $E''(y) = w_0$. In this case $p^i$ is a solution for
all $i$ if $w_0=v_0\cdots v_{d-1}\sim 1$ and for no $i$
otherwise. Similarly if $E''(y) = w_0y^{k_0}$, then $p^i$ is a solution
if and only if $w_0\sim p^{-ik_0}$. In this case the free reduction of
$p^{ik_0}$ is a word of length at least $\abs {ik_0} \abs p_c$ and at
most $\abs {w_0}=K_1$. Consequently $\abs i \abs p_c \le \abs {ik_0}
\abs p_c \le K_1$.

The remaining possibility is that $E''(y) = w_0y^{k_0}\cdots
w_my^{k_m}$ with $m\ge 2$, all $k_j\ne 0$ and no $w_j$ commuting with
$p$. No $w_j$ conjugates $p$ to $p\inv$ either, as $p$ and $p\inv$ are
not conjugate in the free group $F$.  If $p^i$ is a solution, then by
Lemma~\ref{cyclic-reduction} (with $E''(y)$ in place of $E(x)$) some
$p^{ik_j}$ must be a pseudosolution of
$p^{ik_{j-1}}w_{j}p^{ik_j}w_{j+1}p^{ik_{j+1}}$. Lemma~\ref{powers2}
now implies that $\abs i \abs p_c < 7(\abs p + \abs {w_j} + \abs p +
\abs {w_{j+1}} + \abs p)\le 21\abs p + 7K_1$. We have proved the
following lemma.

\begin{lemma}\label{power}
If $rp^is$ is a solution to $E(x)=1$ for some $i$ with $\abs i \abs
p_c > 21\abs p + 7(\abs{u_0\cdots u_{d-1}}+ d\abs{rs})$, then  $rp^is$
is a solution for all $i$.
\end{lemma}

Consider a solution $rp^i sq^jt$ to $E(x)=1$. Define $K_2=21\max\{\abs
p, \abs q\}+ 7(\abs{u_0\cdots u_{d-1}}+ d\abs{rst})$. We will show
that either $\abs i\abs p_c$ or $\abs j\abs q_c$ is no larger than
$K_2d$. Thus each parametric word $rp^i sq^jt$ from
Lemma~\ref{reduction} with two parameters may be replaced by a
collection of parametric words with just one parameter, namely
$rp^{i_0} sq^jt$, $rp^isq^{j_0}t$ with with $\abs{i_0} \abs p_c \le
K_2d$ and  $\abs{j_0}\abs q_c \le C_2d$.

Without loss of generality suppose that $i,j\ge 0$, and $p$ and $q$
are not proper powers. In particular the centralizers in the free
group $F$ of $p$ and $q$ are the cyclic subgroups generated by $p$ and
$q$ respectively.

\begin{lemma}\label{nonconjugate} Suppose $p$ is not conjugate to $q$
  or $q\inv$ and $rp^i sq^jt$ is a solution to $E(x)=1$. Then either $\abs
i\abs p_c$ or $\abs j\abs q_c$ is no larger than $K_2=21\max\{\abs
p, \abs q\}+ 7(\abs{u_0\cdots u_{d-1}}+ d\abs{rst})$.
\end{lemma}

\begin{proof}
First suppose that $s\sim 1$ and take the solution to be
$rp^iq^jt$. Write $E(rp^iq^js)=v_0(p^iq^j)^{\eps_1}v_1\cdots v_{d-1}(p^iq^j)^{\eps_d}v_d$. The $v_k$'s are coefficients; call the
$(p^i)^{\eps_k}$'s and $(q^j)^{\eps_k}$'s powers. Consider how a
coefficient $v_k$ might conjugate the power on one side of itself to
the power on the other side. As $p$ is not conjugate to $q$ or
$q\inv$, $v_k$ would either lie in a subword $v_{k-1} p^iq^j v_k
q^{-j}p^{-i}v_{k+1}$ and centralize $q$ or in a subword $v_{k-1}
q^{-j}p^{-i} v_k p^iq^jv_{k+1}$ and centralize $p$. Consequently $v_k$
is freely equal to a nontrivial power (because $E(x)$ is freely
reduced) of $q$ in the first case and a
nontrivial power of $p$ in the second. $W$ is freely equal to the word
obtained by deleting the powers on either side of $v_k$.

Let $W'$ be the word obtained from $W$ by performing all the deletions
discussed in the previous paragraph. Notice that the first and last
powers of $W$ survive and that the new coefficients are either old
coefficients which do not conjugate their adjacent powers into each
other or products $v_kv_{k+1}\cdots v_{k+m}$  of successive
coefficients whose adjacent powers in $W$ have been deleted. In the
latter case the coefficient is an alternating product of
nontrivial powers of $p$ and $q$.

Since $W'\sim 1$, some power is a pseudosolution in a subword of $W'$
consisting to up to three powers and the coefficients between
them. The sum of the length of the coefficients of $W'$ is the same as
that of $W$, namely $\sum\abs{u_k} + d\abs r + d\abs s$. If $\abs
i\abs p_c$ and $\abs j\abs q_c$ exceed the bound given above, then
Lemma~\ref{powers2}
applies and (as $p$ is not conjugate to $q$ or $q\inv$) implies that
some coefficient conjugates one adjacent power to the other. But this
is impossible either because the coefficient is inherited from $W$ or
because the coefficient is an alternating product of
nontrivial powers of $p$ and $q$, and the conjugation would be a
nontrivial relation satisfied by $p$ and $q$, which generate a free
group of rank two.

It remains to reduce to the case $s\sim 1$. Assume $s\not\sim 1$, and rewrite the
solution as $rp^i(sqs\inv)^j(st)$. One of $\abs i\abs p_c$ or $\abs
j\abs {sqs\inv}_c=\abs j\abs q_c$ is at most $21\max\{\abs
p, \abs {sqs\inv} \}+ 7(\abs{u_0\cdots u_{d-1}}+ d\abs{rst})$.
\end{proof}

Finally, consider a solution $rp^i sq^jt$ to $E(x)=1$ with $p$
conjugate to $q$ or $q\inv$. With appropriate changes to $r$,$s$,$t$
and $j$, $rp^isq^jt$ may be rewritten as $rp^isp^jt$ where $p$ is
cyclicly reduced and $s$ does not commute with $p$.

Define $W=E(rp^isp^jt)=v_0(p^isp^j)^{\eps_1}v_1\cdots
v_{d-1}(p^isp^j)^{\eps_d}v_d$, and argue as before. The coefficients
now include the subwords $s^{\pm 1}$ as well as the $v_k$'s. Consider
how a $v_k$ might conjugate the power on one side of itself to
the power on the other side. Since all powers are powers of $p$, $v_k$
would commute with $p$ and hence would itself be freely equal to a
power of $p$. If $\eps_k\eps_{k+1}=-1$, then $v_k\not\sim 1$ and the powers on
either side cancel. However, if
$\eps_k\eps_{k+1}=1$, then the powers do not necessarily cancel but
combine to form a power $p^{\pm(i+j)}$.

Let $W'$ be the word obtained from $W$ by performing all the deletions
and combinations of powers discussed in the previous paragraph. Notice
that the first and last powers of $W$ survive and that the new
subwords between powers surviving from $W$ are either coefficients
from $W$ which do not conjugate their adjacent powers into each other
or alternating products $s^{\pm 1}p^{k_m}v_ms^{\pm
1}p^{k_{m+1}}v_{m+1}\cdots s^{\pm 1}p^{k_n}v_ns^{\pm 1}s$ where the
$v_j$'s which occur are freely equal to powers of $p$. Further if
$p^{k_j}$ occurs between $s$ and $s\inv$, then $k_j=0$ and $v_j$ is
freely equal to a nontrivial power of $p$, while if $p^{k_j}$ occurs
between two $s$'s or two $s\inv$'s, then $k_j=\pm(i+j)$.

There are two possibilities. First if $\abs{(i+j)}\abs p_c \ge K_2$, then we
may consider the subwords $p^{\pm(i+j)}$ to be powers like the $p^{\pm
i}$'s and $p^{\pm j}$'s surviving from $W$ and the subwords between
the powers as coefficients. Lemma~\ref{powers2}
applies and implies that some coefficient conjugates one adjacent
power to the other and hence is a power of $p$. But this
is impossible either because the coefficient is inherited from $W$ or
because the coefficient is an alternating product of
nontrivial powers of $p$ and $s$, and the conjugation would be a
nontrivial relation satisfied by the subgroup generated by $p$ and
$s$, which is free of rank two.

Second if $\abs(i+j)\abs p_c < K_2$, we take just the  the $p^{\pm
i}$'s and $p^{\pm j}$'s from $W$ to be powers. The coefficients are either
inherited from $W$ or alternating products $s^{\pm 1}p^{k_m}v_ms^{\pm
1}p^{k_{m+1}}v_{m+1}\cdots s^{\pm 1}p^{k_n}v_ns^{\pm 1}s$ as above. In
this case $\abs{(i+j)}\abs p=\abs(i+j)\abs p_c \le K_2$, and the total
length of the
coefficients increases to at most $K_2+ (d-1)\abs{(i+j)}\abs p \le
dK_2$. Lemma~\ref{powers2} applies and yields the following lemma.

\begin{lemma}\label{conjugate} Suppose $p$ is conjugate to $q$
  or $q\inv$ and $rp^i sq^jt$ is a solution to $E(x)=1$. Then $\abs
i\abs p_c$ or $\abs j\abs q_c$ is no larger than $dK_2$.
\end{lemma}

\section{The Algorithm}

The algorithm implicit in the preceding analysis may be described as
follows.

\begin{enumerate}
\item The input is an equation $u_0 x^{\eps_0} u_1 x^{\eps_1} \ldots u_{d -
  1}x^{\eps_{d-1}} = 1$ of degree $d\ge 2$ and with freely reduced
  coefficients from a free monoid $\Sigma^*$ over a set $\Sigma$ of
  generators and their inverses for a free group $F$.
\item Let $L$ be the list of words and parametric words and their
  inverses from
  Lemma~\ref{parametric words}. Rewrite the parametric words so that
  they are either ordinary words or have the one of the forms $rp^is$ or
  $rp^isq^jt$ with $p\not\sim 1 \not\sim q$, $p,q$ not proper powers,
  and in the latter case $sqs\inv\not\sim p^{\pm 1}$.
 \item For each ordinary word $w\in L$ test $E(w)\sim 1$ and
  $E(w\inv)\sim 1$. Remove $w$ from $L$.
\item Replace each parametric word $rp^isq^jt$ with words
  $rp^isq^{j_0}t$ and $rp^{i_0}sq^jt$ for all $i_0,j_0$ with
  $\abs{i_0}\abs p_c \le dK_2$ and  $\abs{j_0}\abs q_c \le dK_2$ where
  $K_2$ is as in Lemma~\ref{nonconjugate}.
\item For each word of the form $w=rp^iq$ in $L$, if
  $E(rp^{i_0}s)\sim 1$ where $i_0$ is the least integer greater than
  $(1/\abs p_c)(21\abs p + 7(\abs{u_0\cdots u_{d-1}} + d\abs{rs})$,
  then $x=rp^is$ is a solution for all $i$, otherwise test
  $E(rp^{i_1}s)\sim 1$ for all $\abs{i_1}< i_0$.
\end{enumerate}

We leave it to the reader to check that our preceding analysis implies
the correctness of the above algorithm. To bound the time complexity
let $\abs L$ be the length the list from Step~2 and $M$ the maximum of
$\abs{rpsqt}$ for each entry $rp^isq^jt$. Note that $M$ is also an
upper bound for the length of the coefficients of $E(x)$ and that the
constant $K_2$ from Lemma~\ref{nonconjugate} is $O(dM)$.

Steps 2 and 3 are accomplished in time $O(M\abs L)$, and Step 4 in
time $O(dK_2M\abs L)$.  Let $L'$ be the augmented list from Step
4. $\abs{L'} = O(dK_2\abs L)=O(d^2M\abs L)$, and each entry in $L'$ has the form
$rp^is$ with $\abs{rps}=O(dMK_2)=O(d^2M^2)$.

For each entry there are $O(M + dM + d^2MK_2)=O(d^3M^2)$ tests
performed in Step 5. The time to test the entry $rp^is$ is linear in
the length of $E(rp^is)$, which is $O(d\abs{rp^is}+\abs{u_1\cdots
u_{d-1}})=O(id^3M^2)=O(dK_2d^3M^2)=O(d^5M^3)$. Thus the total time for
Step 5 is $O(\abs{L'}\cdot (d^3M^2)\cdot (d^5M^3))= O((d^2M\abs
L) \cdot (d^3M^2) \cdot (d^5M^3))=O(d^{10}M^6\abs L)$ . Clearly this
estimate bounds the time of the complete algorithm.

Finally let $m$ be the maximum size of a coefficient in $E(x)$. If
follows from Lemma~\ref{parametric words} that $M=O(m)$ and that $\abs
L = O(dm^3)$. Thus the time complexity of our algorithm is
$O(d^{11}m^9)$.

\end{document}